\documentclass[10pt]{article}
\usepackage{amsmath}
\usepackage{amsthm}
\usepackage{cite}
\usepackage{latexsym}
\usepackage{graphicx, epsfig}
\usepackage{amssymb,amsmath}  
\usepackage{multimedia}
\usepackage{subfigure}
\usepackage{mathrsfs}
\usepackage[usenames,dvipsnames]{pstricks}
 \usepackage{epsfig}
 \usepackage{pst-grad} 
\usepackage{pst-plot} 
\usepackage{mathdots}
\def\ms{\par \medskip}

\def\and{\ \mbox{and} \ }
\def\be{\begin{enumerate}}
\def\ee{\end{enumerate}}
\def\bi{\begin{itemize}}
\def\ei{\end{itemize}}
\def\beqn{\begin{eqnarray*}}
\def\eeqn{\end{eqnarray*}}

\def\F{\mathbb F}
\def\r{\mathbb R}

\def\bc{\begin{center}}
\def\ec{\end{center}}    \def\ed{\end{document}}

\def\stab{\text{stab}}

\newtheorem{ex}{Example}[section]

\newtheorem{thm}[ex]{Theorem}

\newtheorem*{ex*}{Example}
\newtheorem{cor}[ex]{Corollary}
\newtheorem{prop}[ex]{Proposition}
\newtheorem{lem}[ex]{Lemma}

\newtheorem*{qstn*}{Open question}

\begin{document}

\title{Decompositions and eigenvectors of Riordan matrices}
\author{Gi-Sang Cheon\thanks{The author was
supported by the National Research Foundation of Korea (NRF) grant
funded by the Korea government (MSIP) (NRF-2019R1A2C1007518, 2016R1A5A1008055).}~$^{a}$, Marshall M. Cohen$^{b}$ and  Nikolaos
Pantelidis$^{c}$ \\
{\footnotesize $^a$ \textit{Department of Mathematics, Sungkyunkwan
University, Suwon $16419$, Rep. of Korea}}\\
{\footnotesize $^{b}$ \textit{Department of Mathematics, Morgan
State University, Baltimore, MD $21251$, USA }}\\
 {\footnotesize $^{c}$ \textit{School of Science, Waterford Institute of Technology, Ireland}}\\
{\footnotesize gscheon@skku.edu, marshall.cohen@morgan.edu,
nikolaospantelidis@gmail.com}}
\date{}\maketitle
\begin{abstract} Riordan matrices are infinite lower triangular matrices determined by a pair of formal power series over the real or complex field. These matrices have been mainly studied as combinatorial objects with
an emphasis placed on the algebraic or combinatorial structure. The present paper
contributes to the linear algebraic discussion with an analysis of Riordan matrices by means of
the interaction of the properties of formal power series with the linear algebra. Specifically, it is shown that if a Riordan matrix $A$ is an $n\times n$ pseudo-involution then the singular values of $A$ must come in reciprocal pairs. Moreover, we give a complete analysis of existence and nonexistence of the eigenvectors of Riordan matrices. This leads to a surprising partition of the group of Riordan matrices into matrices with three different types of eigenvectors. Finally, given a nonzero vector $v$, we investigate the Riordan matrices $A$ that stabilize the vector $v$, {\it i.e.} $Av=v$.
\end{abstract}


\bigskip
\ms
\section{Introduction}

\mbox{\quad} Triangular matrices appear often in matrix theory, applied linear algebra, combinatorics and also as representations of operators
on spaces of formal power series. For instance, Pascal's triangle can be
represented as an infinite lower triangular matrix
$P=[p_{ij}]_{i,j\ge0}$ by putting the triangle of binomial
coefficients, i.e. $p_{ij}={i\choose j}$ into a matrix.
This matrix
representation is influential in linear
algebra \cite{BP1992,Edelman-Strang}. Using the generalized binomial theorem we see that $j^{\rm th}$ column of the matrix $P$ has the formal power series as its generating function:
\begin{eqnarray}\label{e:pascal}
{1\over 1-x}\left({x\over1-x}\right)^j=\sum_{i\ge j}{i\choose j}x^i,\;\;j\ge0.
\end{eqnarray}

This is one way to view a matrix via columns given by means of associated generating
functions.  In 1991, Shapiro, Getu, Woan and Woodson  \cite{SGWW1991} more generally introduced a special matrix group called the
{\it Riordan group} where it has been proved that the Riordan group unifies many themes in
enumeration. Elements of the group are infinite lower triangular matrices called {\it Riordan matrices} or {\it Riordan arrays}. These matrices are analogously defined in terms of column generating functions as those of the Pascal matrix in the expression (\ref{e:pascal}).
\ms
More formally, let $\F$ be the field of real or complex numbers and let $\F[[x]]$ denote the ring of formal power series (f.p.s.) over $\F$.
 A Riordan matrix $A=[a_{ij}]_{i,j\in {\mathbb N}_0}$ over ${\mathbb F}$ is an infinite lower triangular matrix whose entries are determined by
\begin{eqnarray}\label{e:def}
a_{ij}=[x^i]g(x)F(x)^{j}\;\;{\text{or}}\;\;g(x)F(x)^{j}=\sum_{i\ge
j}a_{ij}x^i
\end{eqnarray}
for some $g,F\in \F[[x]]$ such that $g(0)\ne0, F(0)=0, F'(0)\ne0$, where $[x^i]$ is the coefficient extraction operator and ${\mathbb N}_0=\{0,1,\ldots\}$. As is customary, the Riordan matrix is denoted by $A=\big(g(x), F(x)\big)$ or simply $A=\big(g, F\big)$. By definition, the Pascal matrix $P$ is a Riordan matrix given by $\big({1\over 1-x},{x\over 1-x}\big)$.
\vskip.3pc
Using the {\it fundamental property} of Riordan matrices \cite{SGWW1991} asserting that
\begin{eqnarray}\label{e:ftrm}
\big(g(x), F(x)\big)h(x) = g(x)\cdot (h(F(x)),
\end{eqnarray}
it is shown that the set of all Riordan matrices with entries in ${\mathbb F}$ forms a group under the usual matrix multiplication in terms of generating functions,
\begin{eqnarray}\label{product}
(g,F)(h,L)=\big(g\cdot h(F),L(F)\big).
\end{eqnarray}
This group is called the Riordan group over the field ${\mathbb F}$, and it is denoted by ${\cal R}({\mathbb F})$.
\ms
Riordan matrices have been mainly studied as combinatorial objects, while discussion of
 the algebraic structure does appear, for example, in
 \cite{BarryHennessyPandelitis2020,{CheonKim2013},{CheonKS2008},Cohen2019,Cohen2020,Marshall2017,PandelitisThesis,{JL-Nkwanta2013},PW,Shap} and references there in. Recently, Cheon et al. \cite{CLMPS} established an infinite-dimensional Fr\'{e}chet Lie group and the corresponding Lie algebra on the Riordan group ${\cal R}({\mathbb F})$ from the inverse limit approaches of Riordan groups.
\ms
The present paper contributes to the body of work on the algebra of the Riordan group. Some of the proofs of the known results are new and simpler.
Specifically, we investigate properties of Riordan matrices by means of
formal power series from the viewpoint of linear algebra. More precisely, in Section 2 we give a new factorization of a Riordan matrix in terms of almost-Riordan arrays. Moreover, it is shown that if a Riordan matrix $A$ is an $n\times n$ pseudo-involution then the singular values of $A$ must come in reciprocal pairs.
In Section 3 we describe a complete analysis of existence and
nonexistence of the eigenvectors of Riordan matrices. An eigenvector of a Riordan matrix $(g,F)$ of the form ${\bf{h}}_k=(0, \ldots, 0, h_k, h_{k+1}, \ldots )^T$ with $h_k\ne0$ is called  an {\it eigenvector of level $k$}. When $k = 0$, ${\bf{h}}_0$ is called a
{\it primary eigenvector}. It is normalized if $h_k = 1$. Clearly, the corresponding eigenvalue is $\lambda_k = g_0f_1^k$ where $g_0=g(0)$ and $f_1=[x]F$.  A {\it full set of eigenvectors} is defined to be a set of eigenvectors $\{{\bf h}_0, {\bf h}_1, \ldots, {\bf h}_k,  \ldots \}$, in which every possible eigenvector level is achieved.  Note that a full set of eigenvectors is linearly independent, but is not a basis of the vector space of all infinite sequences in ${\mathbb F}$. As a key result, we prove that the Riordan group ${\cal R}({\mathbb F})$ has a partition (see Theorem \ref{Partition}):
\begin{eqnarray}\label{partition}
{\cal R}({\mathbb F})={\cal R}_{\text {full}}\ \sqcup \ {\cal R}_{\text{none}} \ \sqcup\ \bigsqcup_{k = 0}^\infty {\cal R}_k,
\end{eqnarray}
where ${\cal R}_{\text{full}}, \, {\cal R}_{\text{none}} \  \text{and}\   {\cal R}_k$ are infinite families of Riordan matrices, which matrices respectively have full sets of eigenvectors, no eigenvectors, and a one-dimensional set of eigenvectors, all of which are of level $k$. Finally, in Section 4 we analyze the Riordan matrices in the stabilizer Riordan subgroup for which $A{v}= {v}$ for a given nonzero vector $v$.	

\section{Decompositions of a Riordan matrix}

Riordan matrices with entries of nonnegative integers form a special class of lower triangular matrices
of combinatorial interest. In \cite{PW}, Peart and Woodson
studied a class of Riordan matrices $A$ with triple
factorization, $A= PCF$ where $P$ is a Pascal-type matrix, the
second factor $C$ involves the generating function for the Catalan
numbers, and $F$ involves the Fibonacci generating function. This is
a beautiful factorization from the combinatorial point of view. In
this section, we investigate the matrix decompositions of a Riordan matrix from a linear algebra perspective.

In the sequel, we use the notations ${\cal F}_0$ and ${\cal F}_1$ defined by the sets
\begin{eqnarray*}
{\cal F}_0=\left\{\sum_{i=0}^\infty a_ix^i\in {\mathbb F}[[x]]\;|\;a_0\ne0\right\}\;\;{\rm
and}\;\;{\cal F}_1=x{\cal F}_0.
\end{eqnarray*}
As is well known, $({\cal F}_0,\cdot)$ and $({\cal F}_1,\circ)$ are the groups of invertible f.p.s. with respect to multiplication $\cdot$ and composition
$\circ$, respectively.

\ms

 {\bf Decomposition into almost-Riordan matrices}\;  Recently, Barry \cite{Barry2016} introduced the notion of an {\it almost-Riordan array}. It is an infinite lower-triangular matrix denoted by an ordered triple $(a,g,F)$ of power series $a,g\in {\cal F}_0$ and $F\in{\cal F}_1$ whose first column has the generating function $a$, and the remaining columns starting at the (1,1) position coincide with the Riordan matrix $(g,F)$. For example, $(1,g,F)=[1]\oplus(g,F)$ where $\oplus$ denotes the direct sum of two matrices. We now describe a new class of matrix decompositions that differs from the classical ones e.g. $(g,F)=(g,x)(1,F)$, but are similar in spirit.

\begin{thm} A Riordan matrix $(g,F)\in {\cal R}({\mathbb F})$ can be
factorized into almost-Riordan matrices as follows:
\begin{eqnarray}\label{Factor}
(g,F)=(g,{F/x},x)(1,g,F).
\end{eqnarray}
\end{thm}
\begin{proof}  A Riordan matrix $(g,F)$ can be written as
$$
(g,F)=\left(\begin{array}{ccc}
       g_0&\vline&O\\
       \hline
       g_1&\vline &\\
       g_2&\vline&(gF/x,F)\\
       \vdots&\vline&
      \end{array}\right)=(g,gF/x,F)
      $$
      where $g=g_0+g_1x+g_2x^2+\cdots$. Since
      $(gF/x,F)=(F/x,x)(g,F)$ it follows that
      $$
(g,F)=(g,F/x,x)([1]\oplus(g,F))=(g,F/x,x)(1,g,F),
      $$
as required.
\end{proof}

An $n\times n$ Riordan matrix of $A=(g,F)$ is the truncated Riordan matrix of order $n$, denoted by $(g,F)_n$, which is the leading $n\times n$ principal submatrix of $A = (g,F)$.

\begin{cor} Every $n\times n$ Riordan matrix can be expressed as a product of $n$ almost-Riordan matrices.
\end{cor}
\begin{proof} Let $A_n=(g,F)_n$. Then it follows from (\ref{Factor}) that
$$
A_n=\prod_{k=0}^{n-1}\big(I_k\oplus(g,F/x,x)_{n-k}\big),
$$
as required.
\end{proof}

To illustrate the above proof, we consider the case $n=4$:
      \begin{eqnarray*}
A_4=\left(\begin{array}{cccc}
       g_0&0&0&0\\
       g_1&f_1&0&0\\
       g_2&f_2&f_1&0\\
       g_3&f_3&f_2&f_1
      \end{array}\right)\left(\begin{array}{cccc}
       1&0&0&0\\
       0&g_0&0&0\\
       0&g_1&f_1&0\\
       0&g_2&f_2&f_1
      \end{array}\right)\left(\begin{array}{cccc}
       1&0&0&0\\
       0&1&0&0\\
       0&0&g_0&0\\
       0&0&g_1&f_1
      \end{array}\right)\left(\begin{array}{cccc}
       1&0&0&0\\
       0&1&0&0\\
       0&0&1&0\\
       0&0&0&g_0
      \end{array}\right).
\end{eqnarray*}
\ms

 {\bf Singular Value Decomposition}\; Edelman and Strang
\cite{Edelman-Strang} showed that the singular values of the Pascal
matrix of order $n$ must come in reciprocal pairs $\sigma_i$ and $1/\sigma_i$, for $i=1,2,\ldots,n$. One may ask for which Riordan matrix $A\in GL_n({\mathbb R})$,
the singular values of $A$ must come in reciprocal pairs?

Since $(AB)(AB)^T=A(BB^T)A^T$, it follows that $A$ and $AB$ have the same singular values if and only if $BB^T=I_n$, i.e. $B$ is an orthogonal matrix.
If $B$ is a Riordan matrix with $B\ne \pm I_n=\pm(1,x)$, it can be easily shown that $B$ is orthogonal
if and only if $B$ is $(1,-x)_n$ or $(-1,-x)_n$. A Riordan matrix $A$ is called a {\it
pseudo-involution} if $(AM)^2=I$, {\it i.e.} $AM$ is involution where
 $$ M=\pm(1,-x)=\pm{\rm diag}(1,-1,1,-1,\ldots).$$
 For more information about pseudo-involutions, see \cite{CheonKS2008,Ana2017,JL-Nkwanta2013}.
 \begin{thm} \label{ThmTwo}
 Let $A=(g,F)_n$ be an $n\times n$ pseudo-involution. Then the singular values of $A$ must come in reciprocal pairs $\sigma_i$ and
 $1/\sigma_i$, {\it i.e.} $\sigma_1\ge\sigma_2\ge\cdots\ge 1/\sigma_2\ge 1/\sigma_1$, where $\sigma_n=1/\sigma_1$ and $\sigma_{n-1}=1/\sigma_2$.
 \end{thm}
\begin{proof} Since
$(AM)^2=I_n$ and $M^2=I_n$ we have $A^{-1}=MAM^{-1}$ and
$(A^T)^{-1}=MA^TM^{-1}$. Let $S=AA^T$. Then
\begin{eqnarray*}
S^{-1}&=&(A^T)^{-1}A^{-1}=(MA^TM^{-1})(MAM^{-1})=(MA^T)AM^{-1} \\
&=&(MA^T)AA^T(MA^T)^{-1}=(MA^T)S(MA^T)^{-1}.
\end{eqnarray*}
Thus $S$ is similar to $S^{-1}$, which implies that $S$ and $S^{-1}$
have the same eigenvalues
$\lambda_1\ge\lambda_2\ge\cdots\ge\lambda_n$ where
$\lambda_n=1/\lambda_1$, $\lambda_{n-1}=1/\lambda_2$, and so on. Hence the singular values of $A$
must come in reciprocal pairs $\sigma_i$ and
 $1/{\sigma_i}$ where $\sigma_i=\sqrt{\lambda_i}$.
 \end{proof}

\begin{ex} {\rm The Aigner's directed animal matrix \cite{aign},
\begin{equation*}
\left(g,xg\right) \\
= \left(\begin{array}{c c c c c c c }
	1 &  & & & & &  \\
	1 &1  & & & {\rm O}&& \\
	1 & 2& 1 & & & & \\
	2 & 3& 3& 1 & & &   \\
	4 & 6 & 6 & 4 &1& & \\
	9 & 13 & 13 & 10 &5& 1& \\
	\vdots & \vdots & \vdots &\vdots &\vdots &\vdots  &\ddots  \\	
\end{array}\right)
\end{equation*}
is a pseudo-involution where $g=\frac{1+x-\sqrt{1-2x-3x^2}}{2x}$. Its singular values of the first $6\times 6$ matrix are:
$$\sigma_1\simeq 25.976 \geq \sigma_2 \simeq 2.2139 \geq \sigma_3 \simeq 1.2161 \geq \sigma_4 \simeq 0.82230 \geq \sigma_5 \simeq 0.45169 \geq \sigma_6\simeq0.038497,$$
which are reciprocal pairs as the relations $\sigma_6= \frac{1}{\sigma_1},\sigma_5= \frac{1}{\sigma_2}, \sigma_4= \frac{1}{\sigma_3}$ are satisfied.}
\end{ex}

\begin{thm}\label{e:Th2}
Let $A$ be an $n\times n$ pseudo-involution. Then $A$ has a singular value decomposition of the form $A=U\Sigma V^T$ where $U$ is an $n\times n$ orthogonal matrix diagonalizing $AA^T$, $\Sigma={\rm diag}(\sigma_1,\ldots,\sigma_n)$, $\sigma_1\ge\cdots\ge \sigma_n$ for reciprocal pairs $\sigma_i$ and
 $1/\sigma_i$, and
$V=MUP$ for
\begin{eqnarray}\label{back}
M=(1,-x)_n=\left(\begin{array}{cccc}
       1&0&\cdots&0\\
       0&-1&\ddots&\vdots\\
       \vdots&\ddots&\ddots&0\\
       0&\cdots&0&\pm1
      \end{array}\right)_{n\times n}\;{\rm and}\;\;
P=\left(\begin{array}{cccc}
      0&\cdots&0&1\\
      \vdots&\reflectbox{$\ddots$} &1&0\\
       0&\reflectbox{$\ddots$}&\reflectbox{$\ddots$}&\vdots\\
       1&0&\cdots&0
      \end{array}\right)_{n\times n}.
\end{eqnarray}
 \end{thm}
\begin{proof} By the singular value decomposition of $A$, we may assume that $A=U\Sigma V^T$. Since $A$ is a pseudo-involution it follows from Theorem \ref{ThmTwo} that the singular values $\sigma_1,\ldots,\sigma_n$ of $A$ must come in reciprocal pairs $\sigma_i$ and $1/\sigma_i$. Clearly, $U$ and $V$ are $n\times n$ orthogonal matrices that diagonalize $AA^T$ and $A^TA$, respectively. To complete the proof, we show that $V=MUP$ where $M=(1,-x)$ and $P$ is the $n\times n$ backward identity matrix in (\ref{back}). Since $(AM)^2=I_n$ and $M^2=I_n$ it follows from $A=U\Sigma V^T$
that
\begin{eqnarray}\label{e:eq1}
A=MA^{-1}M=M(V\Sigma^{-1}U^T)M=(M V)\Sigma^{-1}(MU)^T.
\end{eqnarray}
Using $\Sigma^{-2}=P\Sigma^{2}P$, we obtain from (\ref{e:eq1}) that
\begin{eqnarray*}
A^TA&=&(MU)\Sigma^{-1}(MV)^T(M
V)\Sigma^{-1}(MU)^T=(MU)\Sigma^{-2}(MU)^T\\
&=&(MUP)\Sigma^{2}(MUP)^T.
\end{eqnarray*}
Since $A=U\Sigma V^T$ and $U^TU=I_n$, we obtain $A^TA=V\Sigma^2V^T$. By taking $V=MUP$ we thus complete the proof.\end{proof}
\section{Eigenvectors of a Riordan matrix}\label{sec3}

In this section, we study existence of eigenvectors from the viewpoint of the formal power series which generate a Riordan matrix.
 We prove that there is a surprising partition of ${\cal R}({\F})$ given in (\ref{partition})  into three types of Riordan matrices, according to the sets of eigenvectors of these matrices. We also give algorithms for the construction of each type of matrix and give combinatorial criteria for recognizing the type of a given matrix.
\ms

In the sequel, let $A=(g,F)$ be a Riordan matrix in ${\cal R}({\mathbb F})$ where $g=\sum_{n\ge0}g_nx^n\in{\cal F}_0$ and $F=\sum_{n\ge1}f_nx^n\in{\cal F}_1$.
   By the fundamental property in (\ref{e:ftrm}), ${\bf{h}}=(h_0, h_1,\ldots\ )^T$ is an eigenvector of $A$ associated with eigenvalue
$\lambda$  if and only if
\begin{equation} \label{FundThrm}
g(x)h\big(F(x)\big) = \lambda h(x),
\end{equation}
where $h(x)=h_0+h_1x+\cdots$.

 Formal power series $F$ with $F(0) = 0$ play a key role in the section. If $F$ and $ H$ are conjugate in the group ${\cal F}_1$ with operation of composition, we write $F\large \sim H$, in which case there exists $\theta\in{\cal F}_1$ such that $(\theta\circ F\circ\overline{\theta})(x) = H(x)$.
  By $F^{(n)}$, we mean the $n$-fold composition of $F$. We say that $F$ has finite {\it compositional order} $n$, if $n$ is the smallest positive integer such that $F^{(n)}(x)=x$, and that $F$ has infinite compositional order if no such $n$ exists. The finite {\it multiplicative order} of $a\in{\mathbb F}$ is the smallest positive integer $n$ such that $a^n=1$, and if no such $n$ exists, $a$ is said to have infinite multiplicative order.

\begin{lem}\label{lem3.1}{\rm\cite{Cohen2018a,Scheinberg1970}}\; Let $F = f_1 x + f_2x^2 + \cdots\in{\cal F}_1$. Then exactly one of the following  conditions holds:
\begin{itemize}
\item [{\rm(a)}]  $F$ has finite compositional order $n$, in which case $f_1$ has finite multiplicative order $n$ and $F\sim f_1x$.
 \item [{\rm(b)}] $f_1$ has infinite multiplicative order, in which case $F$ has infinite compositional order and $F\sim f_1x$. Indeed there exists a unique series $\theta=\sum_{i=1}^\infty \theta_ix^i$ with $\theta_1=1$ such that $(\theta\circ F\circ\bar{\theta})(x)=f_1x$.
\item [{\rm(c)}] $F$ has infinite compositional order and $f_1$ has finite multiplicative order $n$. Such series $F$ is called a {\rm hybrid series}. This occurs if and only if $F(x)\nsim f_1 x$. Rather, one of the following holds:
\begin{itemize}
\item [{\rm(i)}] if $f_1 = 1$ and $F = x + f_kx^k + \cdots$, $f_k\ne 0$ then there exists unique $b \in \F$ such that
$F \ {\large \sim} \ x + x^k + bx^{2k - 1}$;
\item [{\rm (ii)}] if $f_1^n=1$ ($n$ smallest $\ge2$) then there exist a unique integer $k \geq 2$ and unique $b, c \in \F$  such that
$F \ {\large \sim}\  f_1x + bx^k + cx^{2k - 1}$, where $b\ne 0$ and $k \equiv 1$ (mod $n$).
\end{itemize}
\end{itemize}
\end{lem}

If $f_1$ is a primitive $n^{\rm th}$ root of unity, there exist uncountably many series $F\in {\cal F}_1$ of finite compositional order and also infinitely many hybrid series $F \in {\cal F}_1$. These series can be constructed by setting $F$ as an arbitrary conjugate of the appropriate canonical form given in Lemma \ref{lem3.1}.
 Alternatively, let $F= f_1x + \sum_{k=2}^\infty f_kx^k$ where $f_k$ is arbitrarily chosen for $k\not\equiv 1 \ \text{(mod $n$)}$ and where those $k \equiv 1 \ \text{(mod $n$)}$ do or do not (according to whether $F(x)$ is of finite order  or is hybrid) satisfy certain  required equations  (see \cite{Brewer,Cohen2018a}).
 \ms
 It is convenient in exhibiting Riordan matrices with no eigenvectors or only a one-dimensional set of eigenvectors, to have the most easily recognizable hybrid series, given by the following lemma.

\begin{lem}\label{RecogHybrid} {\rm (Hybrid series)} Suppose that $F \in {\cal F}_1$ is of the form
\[F = f_1x + f_sx^s + f_{s + 1}x^{s + 1} + \cdots, \ (s > 1, \, f_s\ne 0),\]
where $f_1$ has finite multiplicative order $n \geq 1$.\ Then $F$ is a hybrid series if one of the following holds:
\be
\item [{\rm(a)}] $F$ is a polynomial: $F = f_1x + f_sx^s +\cdots + f_qx^q$, $(1 < s \leq q$, $f_s, f_q \ne 0)$.
\item [{\rm(b)}]  $s \equiv 1\ \text{(mod n)}$.
\ee
\end{lem}
\begin{proof} We must prove that $F^{(n)}(x)\ne x$ for all $n\in\mathbb N$. In situation (a) this is true because we see inductively that for any $n\in\mathbb N$,
\[F^{(n)}(x) = f_1^n x + \cdots + f_q^{q^{n - 1}}x^{q^n} \ne x.\]
For (b), it is known (\cite[Lemma 2.4]{Cohen2018a}, \cite[Prop. 2.3.3]{Brewer}) that if $F^{(n)}(x)= x$ and $s \equiv 1$ (mod $n$) for some $n\in\mathbb N$, then $f_s = 0$. Therefore, under our hypothesis $F$ cannot have finite order $n$.
\end{proof}

We also need the following lemma called Roots Theorem in ${\F}[[x]]$.
\begin{lem}\label{roots} {\rm (Multiplicative roots of formal power series)}
\begin{itemize}
\item[{\rm (a)}] {\rm Niven's Theorem} {\rm(\cite[Thm.3]{Niven1969})}\; Suppose that $A = 1 + a_1x + a_2x^2 + \cdots$. Then there exists
a unique series of the form $B = 1 + b_1x + b_2x^2 + \cdots$, such that $B^n = A.$ We denote $B= A^\frac{1}{n}$.
\item[{\rm (b)}] {\rm Extension of Niven's Theorem} {\rm(\cite{Cohen2019, XXD})}\; Suppose that for $k\ge0$,
$$
C=\sum_{n\ge k}c_nx^n=c_kx^k\left(1+\sum_{n\ge1}{c_{k+n}\over
c_k}x^n\right),\;c_k\ne0.
$$
Then, for each $b_1\ne0$, there exists a unique series of the form $B=\sum_{n\ge1}b_nx^n$ such that $B^k= C$ if and only if $\ b_1^k = c_k$ and  $B
= b_1x \widehat{C}^\frac{1}{k}$, where $\widehat{C}=1+\sum_{n\ge1}{c_{k+n}\over
c_k}x^n$.
\end{itemize}
\end{lem}

\subsection{Riordan matrices of  ${\cal R}_{\text{full}}$  and diagonalization}
Let ${\cal L}(\F)$ be the group of all invertible lower triangular matrices over $\F$, which has ${\cal R}(\F)$ as a subgroup. Our opening  examples of matrices with full sets of eigenvectors are given by the following Proposition.

\begin{prop}\label{uniLTeig} If $A =[a_{ij}]_{i, j \in{\mathbb N}_0} \in {\cal L}(\F)$ has distinct diagonal elements then for each integer $k \geq 0$ there exists a unique eigenvector of level $k$ of the form ${\bf h}_k = (0,\,  \ldots,\, 1, \, h_{k +1}, \ldots )^T$ associated with eigenvalue $\lambda = a_{kk}$.
\end{prop}
\begin{proof} Fix $k \geq 0$. Let $\lambda = a_{kk}$ and $A_{(n)}$ the $n^{\rm th}$ row of the matrix $A$. We consider a variable vector ${\bf h} = (0, \ldots, 0, 1, h_{k + 1}, \ldots)$ of level $k$. Note that the dot product $\langle A_{(n)},{\bf h}\rangle = 0 = \lambda h_n$ for $n < k$, and $\langle A_{(k)},{\bf h}\rangle = a_{kk}\cdot 1 = \lambda h_k$.
We solve the equation $A{\bf{h}}=\lambda{\bf{h}}$ inductively for $h_n$ with $n \geq k + 1$. Thus we seek solutions to
\[0 = \langle A_{(n)},{\bf h}\rangle - \lambda h_n =
 {\textstyle \sum_{j = 0}^{n - 1} a_{n j} h_j}\ + (a_{nn} - a_{kk})h_n.\]
By hypothesis, $a_{nn} \neq a_{kk}$. Thus we may solve uniquely for $h_n$,  proving the theorem.
\end{proof}
\begin{cor}\label{UniEig}\mbox{}
Let $A=(g, F)$ be a Riordan matrix in ${\cal R}(\F)$. If $f_1$ has infinite multiplicative order in $\F\setminus \{0\}$ then for every $k\ge 0$, $A$ has a
unique eigenvector of the form\
 ${\bf{h}}_k = \big(0,\ldots, 0,1,\  h_{k+1},\ldots \big)^T.$
\end{cor}
\begin{proof} If $i \neq j$ then $a_{ii} = g_0f_1^i \neq g_0f_1^j = a_{jj}$. We may apply the Proposition \ref{uniLTeig}.
\end{proof}

By definition, $A \in {\cal L}(\F)$ is {\it diagonalizable} if and only if there exists a diagonal matrix $D \in  {\cal L}(\F)$ such that $A$ is conjugate to $D$ in the group ${\cal L}(\F)$, {\it i.e.} there exists a diagonalizing matrix $X \in {\cal L}(\F)$ such that $X^{-1}AX = D$.
 The classical connection of diagonalizability to the existence of a full set of eigenvectors corresponding to eigenvectors for $A\in {\cal L}(\mathbb F)$ is given by the Lemma below.
\begin{lem} \label{LTeigencolumns} $A \in {\cal L}(\F)$ is diagonalizable if and only if $A$ has a full set of eigenvectors.
\end{lem}
\begin{proof} Let $X\in {\cal L}(\F)$ be a diagonalizing matrix such that $X^{-1}AX = D$. Clearly, the columns of $X$ form a full set of eigenvectors for $A$, and the diagonal elements of $D$ are the corresponding eigenvalues of these eigenvectors of $A$. Since a full set of eigenvectors of $A$ is linearly independent it follows the converse.
    \end{proof}

To recognize whether a diagonalizing matrix of a diagonalizable Riordan matrix can in fact be taken to be a Riordan matrix, we will  need the following Lemma. The proof is straightforward.
\begin{lem} \label{FunctEq} Let $A =(g,F)$ and $X =(h, \theta)$ be Riordan matrices. Then $X^{-1}AX =(g_0,f_1x)$ if and only if $X$ gives a solution to the functional equations:
\begin{eqnarray}\label{solutions}
g_0h(x)=g(x)h(F(x))\quad{and} \quad f_1\theta(x)=\theta(F(x)).
\end{eqnarray}
\end{lem}
Thus, $A$ is diagonalizable in ${\cal R}(\F)$ if and only if $A$ has a primary eigenvector $h$ and there exists $\theta \in {\cal F}_1$ such that $(\theta\circ F\circ\overline{\theta})(x) = f_1x$.
In the light of Lemma \ref{FunctEq}, our key tool in forming a unified theory is the following
\begin{lem} \label{2Vec} {\rm(Two Vector Lemma)}\; If $A = (g, F)$ has two linearly independent eigenvectors then  there exists $\theta \in {\cal F}_1$ such that $(\theta\circ F\circ\overline{\theta})(x) = f_1x$.
\end{lem}
\begin{proof}

Let ${\bf u} = (0, \ldots, 0, u_\ell, u_{\ell + 1}, \ldots )$ and ${\bf v} = (0, \ldots, 0, v_k, v_{k + 1}, \ldots )$ be linearly independent eigenvectors of $A$.
We may assume that  $u_\ell = 1$ and $v_k = 1$.  If $k = \ell$ then  ${\bf u}$ and $\bf {v}$ have the same eigenvalue $g_0f_1^k$, so that ${\bf u}$ and $\bf {v - u}$ are eigenvectors of different levels. Thus we may assume that $k > \ell \geq 0$.

To find $\theta \in {\cal F}_0$ which conjugates $F$ to $f_1x$, we let
$u(x) \ \text{and}\  v(x)$ be the generating functions of the eigenvectors $\bf u$ and $\bf v$. We have the system
\begin{eqnarray*}
 g(x)\cdot u\big(F(x)\big) &=& g_0f_1^\ell\cdot u(x),\\
 g(x)\cdot v\big(F(x)\big) &=& g_0f_1^k\cdot v(x).
\end{eqnarray*}
By substituting $a(x) := \frac{v(x)}{u(x)}$ we obtain $a\big(F\big)={v(F)\over u(F)}=f_1^{k-\ell}a(x)$. Using Niven's Theorem, we write
\begin{eqnarray}\label{ax}
a(x)=x^{k - \ell}(1 + a_1x  + a_2x^2 + \cdots)={\theta}(x)^{k - \ell},
\end{eqnarray}
where ${\theta}(x)=x(1 + a_1x  + a_2x^2 + \cdots)^{1/(k - \ell)}$. By substituting $x=F$ in this equation we obtain
$$
{\theta}(F)^{k - \ell}=a(F)=f_1^{k-\ell}a(x)= f_1^{k - \ell}{\theta}(x)^{k - \ell}.
$$
By the extension of Niven's Theorem in Lemma \ref{roots}, we have ${\theta}\big(F(x)\big)=f_1{\theta}(x)$, which implies
$$
({\theta}\circ F\circ\overline{\theta})(x)={\theta}(F(\overline{\theta}(x)))=f_1{\theta}(\overline{\theta}(x))=f_1x,
$$
as desired.\end{proof}
\begin{lem}\label{2Vec-primary} If $A = (g, F)$ has an eigenvector of level $k > 0$ and if there exists $\theta \in {\cal F}_1$, such that $(\theta\circ F\circ\overline{\theta})(x) = f_1x$, then $A$ has also a primary eigenvector.
\end{lem}
\begin{proof}
Let $v(x) = x^k(1 + v_1x + \cdots)$ give an eigenvector of level $k$. Let $\theta(x) = \theta_1x + \theta_2x^2 + \cdots$. Consider the eigenvector ${\bf h}$ given by
\begin{eqnarray}\label{hx}
h(x) := v(x)\theta(x)^{-k} = \theta_1^{-k}  + h_1x + \cdots
\end{eqnarray}

  Since $v(x)$ generates an eigenvector of level $k$, it follows that $(g,F)v(x) = g(x)v\big(F\big) = g_0f_1^kv(x)$.  
 Also, $\theta(x)$ satisfies ${\theta}(F(x)) = f_1{\theta}(x)$. Thus we have
 \begin{eqnarray*}
(g,F)h(x)  & = & (g,F)v(x){\theta}(x)^{-k}=g(x)v(F){\theta}(F)^{-k}=g_0f_1^kv(x)(f_1{\theta}(x))^{-k} \\
 & = & g_0f_1^{-k + k}v(x){\theta}(x)^{-k}=g_0h(x),
\end{eqnarray*} which proves that $h(x)$ generates a primary eigenvector of $A$.
\end{proof}

From the results above,  we can now paint the total picture with a set of equivalent criteria for recognizing when $(g,F)\in {\cal R}_{\text{full}}$. One is that $(g,F)$ is conjugate in ${\cal R}({\mathbb F})$ to $(g_0, f_1x)$, which solves the diagonalizability problem for Riordan matrices.

\begin{thm}\label{FullEquiv} {\rm (Riordan matrices in ${\cal R}_{\text{full}}$)} Let $A=(g,F)\in{\cal R}(\mathbb F)$. Then the following statements are equivalent:
\be
\item [{\rm(a)}] $A$ has a full set of eigenvectors, i.e. $A \in {\cal R}_{\rm full}$.
\item [{\rm(b)}] $A$ has two linearly independent eigenvectors.
\item[{\rm(c)}] $A$ has an eigenvector of level $k>0$ and there exists $\theta \in {\cal F}_1$ such that $(\theta\circ F\circ\overline{\theta})(x) = f_1x$.
\item [{\rm(d)}] $A$ has a primary eigenvector and there exists  $\theta\in {\cal F}_1$ such that $(\theta\circ F\circ\overline{\theta})(x)= f_1x$.
\item[{\rm(e)}] $A$ is conjugate to the diagonal matrix $(g_0,f_1x)$ by a Riordan matrix $(h, \theta)\in{\cal R}(\mathbb F)$ such that
\begin{eqnarray}\label{R-conjugate}
 (h, \theta)^{-1}(g, F)(h, \theta)=(g_0, f_1x).
\end{eqnarray}
\item[{\rm(f)}] $A$ has an eigenvector of level $k$ for some $k \geq 0$ and either, $F$ has finite compositional order in ${\cal F}_1$, or $f_1$ has infinite multiplicative order in ${\mathbb F}$.

\item[{\rm(g)}] Let $\widehat{g} := \frac{1}{g_0}g(x)\in{\cal F}_0$. Either $f_1$ has infinite multiplicative order in $\F\setminus \{0\}$ or  $\big(\widehat{g},F\big)$ is a finite order element of the Riordan group ${\cal R}(\mathbb F)$ of order equal to the order of $f_1$. \ee
\end{thm}
\begin{proof} Clearly (a) $\implies$ (b). The implications, (b)$\implies$(c)$\implies$(d)$\implies$(e)$\implies$(a), follow in turn directly from Lemmas \ref{2Vec}, \ref{2Vec-primary}, \ref{FunctEq}, \ref{LTeigencolumns}. Moreover, (c)$\iff$(f) follows from the classification of elements of ${\cal F}_1$ under conjugation; see Lemma \ref{lem3.1}.
\ms
\noindent Finally, we prove that (e)$\implies$(g) and (g)$\implies$(a):\\
(e)$\implies$(g): Suppose that $f_1$ does not have infinite order, but in fact has finite order $n$. Notice that $(g,F) = g_0\big(\widehat{g}, F\big)$ and $(g_0, f_1x) = g_0(1, f_1x)$.  Since the matrix $(h, \theta)^{-1}$ gives a linear transformation,  our assumption (e) yields
$$
g_0(h, \theta)^{-1}\,\big(\widehat{g}, F\big)\,(h, \theta) = (h, \theta)^{-1}\,\big(g, F\big)\,(h, \theta)=(g_0, f_1x) = g_0(1, f_1x),
$$
which implies $(h, \theta)^{-1}\,\big(\widehat{g}, F\big)\,(h, \theta) =(1, f_1x).$
The order of $(1, f_1x)$ in ${\cal R}(\mathbb F)$ is finite, since it equals the order $n$ of $f_1$ in $\F\setminus \{0\}$. But conjugate elements in a group have the same order.  Thus $(\widehat{g}, F)$ is of finite order in ${\cal R}(\mathbb F)$.\\
(g)$\implies$(a): If $f_1$ has infinite order, then $(g,F) \in {\cal R}_{\rm full}$ by Corollary \ref{UniEig}. Now let $\big(\widehat{g}, F\big)$ be of finite order. Then Theorem 3 of \cite{Cohen2020} exhibits a full set of eigenvectors for $\big(\widehat{g}, F\big)$. But  $\bf h$ is an eigenvector of level $k$ for $\big(\widehat{g}, F\big)$  if and only if $g_0{\bf h}$ is an eigenvector of level $k$ for $g_0\big(\widehat{g}, F\big) =(g,F)$. Therefore $(g,F)$ has a full set of eigenvectors.
\end{proof}

\begin{ex} {\rm Let
$$
A=\left(\frac{1+x}{1 - x},\, -x\right)=\left(\begin{array}{ccccc}
       1&&&&\\
       2&-1&&{\rm O}&\\
       2&-2&1&&\\
       2&-2&2&-1&\\
       \vdots&\vdots&\vdots&\vdots&\ddots
      \end{array}\right).
$$

Clearly, ${\bf h}_0=(1,1,0,\ldots)^T$ is a primary eigenvector of $A$ associated with $\lambda=1$. Moreover, $\theta = x$ conjugates $-x = F$ to $-x = f_1x$. Thus it follows from the equivalence of (d) and (e) in Theorem \ref{FullEquiv} that,
 taking $(h, \theta) = (1 + x, x)$ with  $(h, \theta)^{-1} = \left(\frac{1}{1 + x}, x\right)$ we may apply (e) to get:}
  {\small \begin{eqnarray*}
\left(\begin{array}{cccc}
    \  \ 1&\ \ 0&\ \ 0&0\\
        -1&\ \ 1&\ \ 0&0\\
   \  \  1&-1&\ \ 1&0\\
        -1&\ \ 1&-1&1\\
       &\cdots&\cdots&
      \end{array}\right)\left(\begin{array}{cccc}
       1&0&0&0\\
       2&-1&0&0\\
       2&-2&1&0\\
       2&-2&2&-1\\
        &\cdots&\cdots&
      \end{array}\right)\left(\begin{array}{cccc}
       1&0&0&0\\
      1&1&0&0\\
       0&1&1&0\\
       0&0&1&1\\
        &\cdots&\cdots&
      \end{array}\right)=\left(\begin{array}{cccc}
       1&0&0&0\\
       0&-1&0&0\\
       0&0&1&0\\
       0&0&0&-1\\
        &\cdots&\cdots&
      \end{array}\right).
\end{eqnarray*}}
\end{ex}

A fundamental question is that if a Riordan matrix $A \in {\cal R}(\F)$ is diagonalizable in ${\cal L}(\F)$, then is it diagonalizable in ${\cal R}(\F)$? This will be answered affirmatively in the following theorem.
\begin{thm}\label{LT-diag} Let $A =(g,F)\in {\cal R}(\mathbb F)$. Then $A$ is diagonalizable in ${\cal L}(\F)$ if and only if $A$ is diagonalizable in ${\cal R}(\F)$.
\end{thm}
\begin{proof} The sufficiency is immediate. For necessity, assume that there exists a $X \in {\cal L}(\mathbb F)$ such that $X^{-1}AX=(g_0,f_1x)$.
Then the columns of $X$ form a full set of eigenvectors for $A$ and the diagonal elements of $(g_0,f_1x)$ are the corresponding eigenvalues of these eigenvectors of $A$. Thus Theorem \ref{FullEquiv} implies that $A$ is diagonalizable in ${\cal R}(\F)$.
\end{proof}

Now, we prove the partition (\ref{partition}) of the Riordan group ${\cal R}(\mathbb F)$ that we stated in Section 1.
\begin{thm} \label{Partition} The Riordan group ${\cal R}(\mathbb F)$ can be partitioned into matrices with three
different types of eigenvectors given by
$${\cal R}(\mathbb F)= {\cal R}_{\rm full}\, \sqcup \, \bigsqcup_{k = 0}^\infty {\cal R}_k \, \sqcup\,  {\cal R} _{\rm none}.$$
\end{thm}
\begin{proof} By definition, the sets ${\cal R}_{\rm full}, \big({\cal R}_k\big)_{k \geq 0}\ \text{and}\ {\cal R} _{\rm none}$ are pairwise disjoint. It is also easily shown that these sets are nonempty. Let $A=(g, F) \in {\cal R}(\mathbb F)$. Then $A$ is an element of ${\cal R}_{\rm full}$ or $A$ is not in ${\cal R}_{\rm full}$. If $A\notin {\cal R}_{\rm full}$ then by Theorem \ref{FullEquiv}, part (b), either no eigenvectors exist or all those which exist are multiples of a single vector. If the level of this vector is $k$ for some integer $k$, then $A \in {\cal R}_k$.   Thus $A \in \bigsqcup_{k = 0}^\infty {\cal R}_k \, \sqcup\,  {\cal R} _{\rm none}$, so that ${\cal R}(\mathbb F)\subset {\cal R}_{\rm full}\, \sqcup \, \bigsqcup_{k = 0}^\infty {\cal R}_k \, \sqcup\,  {\cal R} _{\rm none}$. The reverse is clear. Thus we have a partition of ${\cal R}(\F)$.
\end{proof}

\subsection{Construction of Riordan matrices in ${\cal R}_{\text{full}}, {\cal R}_k, \ \text{and}\ {\cal R} _{\rm none} $}

From Theorem \ref{FullEquiv} we get complete prescriptions for constructing elements $(g, F)$ of ${\cal R}_{\text{full}}$, ${\cal R}_k$ and ${\cal R}_{\text{none}}$. In making these constructions when $f_1$ is a primitive $n^{\rm{th}}$ root of unity, $F\in{\cal F}_1$ will either be of finite compositional order or will be hybrid series according to the criteria given in Lemmas \ref{lem3.1} and \ref{RecogHybrid}. Having chosen $F$, it follows from (\ref{FundThrm}) that $g$ will be determined by any eigenvector $h(x)$.

\begin{thm}\label{exRfull} Every element $(g, F) \in {\cal R}_{\text{full}}$ can be constructed as follows:
\be
\item [{\rm(i)}] Choose $g_0 \ne 0$.
\item [{\rm(ii)}] Choose a conjugate $F$ in ${\cal F}_1$ of  some $f_1x\ (f_1\ne0)$.
\item [{\rm(iii)}] Choose a series $h = h_kx^k + \cdots\, (h_k \ne 0)$ of some level $k \geq 0$.
\item  [{\rm(iv)}] Set $g = g_0f_1^k\cdot\frac{h(x)}{h\big(F(x)\big)}$.
   \ee
 \end{thm}
\begin{proof} Applying Theorem \ref{FullEquiv} (c), the result follows.
\end{proof}

In the following theorem, we construct, for each $k \geq 0$ the  class ${\cal R}_k$ of Riordan matrices which have eigenvectors of level $k$ and of no other level.

  \begin{thm}  \label{R_kElement} Every element $(g, F) \in {\cal R}_k$ can be constructed as follows:
  \be
 \item [{\rm(i)}]  Choose $g_0\ne0$.
\item [{\rm(ii)}] Choose a hybrid series $F$ (Lemma 3.1, Lemma \ref{RecogHybrid}).
\item [{\rm(iii)}] Choose a series $h = h_kx^k + \cdots\ ( h_k\ne 0)$ of some  level $k \geq 0$.
\item [{\rm(iv)}] Set $g = g_0f_1^k\cdot\frac{h(x)}{h\big(F(x)\big)}$
\ee
Moreover, if $(g, F) \in {\cal R}_k$ then there is a unique $h$ of the form
 $h = x^k + h_{k + 1}x^{k + 1} + \cdots$ such that $g = g_0f_1^k\frac{h(x)}{h\big(F(x)\big)}$ .
 \end{thm}
\begin{proof} Statements (iii) and (iv) are equivalent to the statement that $(g, F)$ has an eigenvector of level $k$. Given this, (ii) implies that $(g, F)$ has eigenvectors of no other level, by the equivalence of (b) and (d) of Theorem \ref{FullEquiv}.  The uniqueness follows from Theorem \ref{FullEquiv} (b) that two independent eigenvectors would imply the existence of a full set of eigenvectors, contradicting $(g, F)  \in {\cal R}_k$. \end{proof}

\begin{ex} \label{R_kExamples} {\rm The simplest examples of elements of ${\cal R}_k$ given by Theorem \ref{R_kElement} are constructed by letting $g_0 = 1, \, F = \pm x + x^2$ and $h = x^k$.  This gives the examples respectively
\[ A_k = \left(\frac{1}{(1 + x)^k},\ x + x^2\right) \in  {\cal R}_k\ \  \text{and} \ \ B_k = \left(\frac{1}{(1 - x)^k},\ -x + x^2\right) \in  {\cal R}_k\]
For instance, $B_1 =  \big(\frac{1}{1 - x}, - x + x^2\big)$ has eigenvector $(0,1,0,\ldots)^T$ of level one associated with $\lambda=-1$ and no eigenvectors of any other level.}
\end{ex}

In the following theorem we show that, from Theorem \ref{FullEquiv}, we can construct surprising elements  $(g,F)\in {\cal R}_{\text{none}}$, because $(g, F)$ has no eigenvectors even though $F$ is of finite compositional order and is thus conjugate to $f_1x$ in ${\cal F}_1$.

\begin{thm}\label{NoneDespite} Suppose that $F\in {\cal F}_1$ has finite compositional order $n$ and that $g =1 + g_rx^r + g_{r + 1}x^{r + 1} + \cdots\in {\cal F}_0$ with $r > 0$,  $g_r\ne0$.
If $f_1^r = 1$, i.e. $r \equiv 0$ (mod $n$), then $(g,F)\in {\cal R}_{\text{none}}$.
\end{thm}
\begin{proof} If $(g, F)$ had any eigenvectors whatsoever, it would follow from the equivalence of (d) and (e) in Theorem \ref{FullEquiv} that $(g, F)^n = (1, x)$.  Thus, looking at the first coordinate, we have
\[g\cdot g\big(F\big) \cdots \,\cdot g\big(F^{(n - 1)}\big) = 1.\]
Noticing, by  induction, that $F^{(j)}(x) = f_1^jx + (\text{higher powers})$, we obtain
$$
g\big(F^{(j)}\big) = 1 +\sum_{k\ge0} g_{r+k}\big(f_1^jx  + (\text{higher powers})\big)^{r+k}.
$$
Therefore
\beqn 
0 & = &  [x^r]\,1 = [x^r]\ \Big(g\cdot g\big(F\big) \cdots \cdot g\big(F^{(n - 1)}\big) \Big)\\
& = & [x^r] \Big((1 + g_rx^r)(1 + g_rf_1^rx^r)\cdots (1 + g_r(f_1^{n - 1})^rx^r\Big)\\
& = & [x^r]\ g_r\,\Big( 1 + (f_1^r) + \cdots+ (f_1^r)^{n - 1}\Big)x^r,\ \text{where}\ f_1^r = 1,\\
& =& n g_r \neq 0.
\eeqn
This contradiction implies that $(g, F)$ has no eigenvectors.
\end{proof}

\begin{ex} {\rm Let $A= \Big(1 + g_2x^2 + g_3x^3 + \cdots,\, -x \Big)$. If $g_2 \ne0$ then $A$ has no eigenvectors.}
\end{ex}
  From Theorem \ref{FullEquiv} we see that $(g, F) \in
 {\cal R} _{\rm none} \sqcup \, \bigsqcup_{k = 0}^\infty {\cal R}_k$ \ if and only if \ $(g, F)$ does not have two linearly independent eigenvectors.  For hybrid series $F$,  statement (d) of Theorem \ref{FullEquiv} is false. Thus hybrid series play the following role:
\begin{thm}\label{hybridEig} If $F$ is a hybrid series then
$(g, F) \in {\cal R}_{\rm none}$ or $(g, F) \in {\cal R}_{\rm k}$ for some unique integer $k \geq 0$.
\end{thm}

The converse of Theorem \ref{hybridEig} is false. We have shown in Theorem \ref{NoneDespite} that there exist  $(g, F) \in  {\cal R}_{\text{none}}$ such that $F$ is not a hybrid series, but rather a series of finite compositional order.
In general, it is difficult to recognize, for a given element $(g, F)$ which does not have two independent eigenvectors, whether or not
$(g, F)$  has an eigenvector at all. Thus we have a question: given $(g, F)$ can we use quick numerical computation to recognize the type of $(g, F)$?
The next subsection gives results on this question.

\subsection{Recognizing the eigenvector type of a given element $(g, F)\in {\cal R}(\mathbb F)$} \label{Recogn}

We begin by applying the results above to record the answers in the most easily recognizable examples -- those of the forms $\big(g, F\big)= (g_0, F)$ or $\big(g, \, f_1x\big)$.

\begin{thm}\label{Recognizable} Let $(g_0, F)\in{\cal R}(\mathbb F)$ where $F = f_1x + f_2x^2 + \cdots\in {\cal F}_1$.
\begin{itemize}
\item [{\rm (a)}] If $f_1$ is of infinite multiplicative order then $(g_0, F) \in {\cal R}_{\rm full}$.
\item [{\rm (b)}] Assume that $f_1$ is of finite multiplicative order. If $F$ is of finite compositional order then $(g_0, F) \in {\cal R}_{\rm full}$; and if $F$ is of infinite compositional order then $(g_0, F) \in {\cal R}_0$,
with eigenvectors $(1, 0,\ldots\,)^{\rm T}$ and its multiples.
\end{itemize}
\end{thm}
\begin{proof} The statement (a) follows from Corollary \ref{UniEig}. (b) Let $f_1$ be of finite multiplicative order. If $F$ is of finite compositional order,  by Lemma \ref{lem3.1} (a) we have $F\sim f_1x$. Thus by Theorem \ref{FullEquiv} (e) it follows from $(g_0, F) \sim (g_0, f_1x)$ that $(g_0, F) \in {\cal R}_{\rm full}$. If $F$ is of infinite compositional order then by inspection, $(1, 0,\ldots\,)^{\rm T}$ is an eigenvector of $(g_0, F)$. If there was another eigenvector linearly independent of this, then Theorem \ref{FullEquiv}, parts (b) and (d) would imply that $F$ is not hybrid. Thus $(g_0, F) \in {\cal R}_0$.
\end{proof}

\begin{thm}\label{Recognizable-1}  Let $(g,f_1x)\in{\cal R}(\mathbb F)$ where $g=g_0 + g_rx^r + g_{r + 1}x^{r + 1} + \cdots\in {\cal F}_0$ with $g_r\ne0$.
\begin{itemize}
\item [{\rm (a)}] If $f_1$ is of infinite multiplicative order then $\big(g, \,f_1x\big) \in {\cal R}_{\rm full}$.
\item [{\rm (b)}] If $f_1$ is of finite multiplicative order $n$ then either $\big(g, \,f_1x\big)\in  {\cal R}_{\rm full}$ or ${\cal R}_{\rm none}$.
\end{itemize}
\end{thm}
\begin{proof} The statement (a) follows from Corollary \ref{UniEig}. (b) Let $f_1$ be of finite multiplicative order $n$.
Since $F = f_1x$ is not a hybrid series, Theorem \ref{R_kElement} implies $(g, F) \notin {\cal R}_k$. Then, by Theorem \ref{FullEquiv} (f) and (g), $(g, F) \in {\cal R}_{\rm full}$ if and only if $\big(\widehat{g}, f_1x\big)$ has order $n$ where $\widehat{g}= \frac{1}{g_0}g(x)$. This occurs if and only if $(1, x) = \big(\widehat{g}, f_1x\big)^n = \big(\widehat{g}(x)\cdot \widehat{g}\big(f_1x\big) \cdots \widehat{g}\big(f_1^{n-1}x\big),\ f_1^nx\big)$. Thus the statement (b) follows.
\end{proof}

We assume from now on that $(g,F) \in {\cal R}(\mathbb F)$ satisfies
\begin{eqnarray}\label{rs}
\begin{cases}
g=  g_0 + g_rx^r + g_{r + 1}x^{r + 1} + \cdots,\; (r\ge1,\;g_r \neq 0); \\
F= f_1x + f_sx^s + f_{s + 1}x^{s+1} + \cdots,\; (s\ge2,\;f_s \neq 0);\\
f_1\ {\rm has\ finite \ multiplicative \ order}\ n \in{\mathbb N}.
\end{cases}
\end{eqnarray}

The following computational lemma will help us to more precisely recognize whether $(g, F) \in {\cal R}_{\rm none}$ or $(g, F) \in {\cal R}_k$ for some particular $k$.
\begin{lem}\label{recognitionEqns} Let $A=(g, F) \in {\cal R}(\mathbb F)$ where $g, F$ are given in (\ref{rs}). Suppose that $A$ has an eigenvector $\bf h$ of level $k \geq 0$ with generating function \
$h = h_kx^k + h_{k + 1}x^{k + 1} + \cdots.$
\be
\item[{\rm (a)}] If $r < s$ then
\begin{equation}\label{eqnaa}
g_0f_1h_{k + r}(1 - f_1^r)  =  (kg_0f_{r + 1} + g_rf_1)h_k.
\end{equation}
\item[{\rm (b)}] If $r \geq s$ then
\begin{equation}\label{eqnbb}
kh_kf_s = f_1h_{k + s -1}\left(1 - f_1^{s -1}\right).
\end{equation}
\ee
\end{lem}
\begin{proof} First note that the coefficients of $h(F(x))$ can be determined by
\begin{eqnarray}\label{power}
[x^n]h(F(x))=\sum_{j_1+\cdots+j_i=n}h_if_{j_1}\cdots f_{j_i}
\end{eqnarray}
where the sum runs over all positive integer solutions $j_1,\ldots,j_i$ to $j_1+\cdots+j_i=n$ for each $i=1,\ldots,n$.

\noindent(a) Let $r < s$. Since the eigenvalue for {\bf h} is $g_0f_1^k$, we have
$$(g,F)h=gh(F)=g_0f_1^kh.$$
Using (\ref{power}) together with $g_1=\cdots=g_{r-1}=0$ it can be shown that
\begin{eqnarray*}
 g_0f_1^kh_{k + r} & = &  [x^{k + r}] \Big(gh(F)\Big)=\sum_{t=0}^{k+r}[x^t]g\cdot[x^{k+r-t}]h(F)\\
 &=& g_0\big(h_k\cdot kf_1^{k - 1}f_{r + 1}+h_{r + k}f_1^{r + k}\big)+g_rh_kf_1^k.
\end{eqnarray*}
Dividing both sides by $f_1^{k-1}$ yields (\ref{eqnaa}).

\noindent(b) Let $r \geq s$.  Then $g _1 = \cdots = g_{s - 1} = 0$. By a similar method used in (a), we have
\begin{eqnarray*}
g_0f_1^kh_{k + s -1} & = &  [x^{k + s-1}] \Big(g\cdot h(F)\Big)=\sum_{t=0}^{k+s-1}[x^t]g\cdot[x^{k+s-1-t}]h(F)\\
 &=& g_0\big(h_kkf_1^{k - 1}f_s+h_{k + s - 1}f_1^{k + s -1}\big),
\end{eqnarray*}
which gives (\ref{eqnbb}).
\end{proof}

Lemma \ref{recognitionEqns} leads to the following theorem.
\begin{thm}\label{recognitionThm} {\rm(Recognition Theorem)} Let $(g, F) \in {\cal R}(\mathbb F)$ where $g, F$ are given in (\ref{rs}).
 \begin{itemize}
 \item [{\rm (a)}] If $r < s - 1$ and $f_1^r = 1$ then $(g, F) \in {\cal R}_{\rm none}$.
  \item [{\rm (b)}] Let $r = s - 1$ and $f_1^r = 1$. Then the following holds:
  \begin{itemize}
 \item [{\rm (i)}] If there exists an eigenvector of level $k$, then
$\displaystyle{k = -\frac{g_rf_1}{g_0f_s}\ne 0}$, and $(g,F) \in \mathcal{R}_{k}$.
 \item [{\rm (ii)}] If $g_0, g_r, f_1, \ \text {and}\  f_s$ are real numbers with $g_0g_rf_1f_s > 0$, then $(g, F) \in {\cal R}_{\rm none}$.
   \end{itemize}
     \item [{\rm(c)}] Let $r \geq s$ and $f_1^{s - 1} = 1$. If there exists an eigenvector of level $k$, then $k = 0$.
 \end{itemize}
\end{thm}
\begin{proof} (a) If $r < s - 1$ then ${r + 1} < s$, so that $f_{r + 1} = 0$. Thus, if
$f_1^r = 1$, then equation (\ref{eqnaa}) implies $0 = g_rf_1h_k$, contradicting $g_r, f_1, h_k\ne0$. Thus  $(g, F) \in {\cal R}_{\rm none}$.

\noindent (b) Given that $f_1^r = 1$,  (i) follows directly from equation (\ref{eqnaa}).
The result of (i) and the hypothesis of (ii) imply that $k<0$.  This is impossible, so that in the situation of (ii), $(g, F) \in {\cal R}_{\rm none}$.

\noindent (c) If $r \geq s$ and $f_1^{s - 1} = 1$ then equation (\ref{eqnbb}) gives $kh_kf_s = 0$. Thus if there exists an eigenvector, it must be of level $k = 0$, {\it i.e.} a primary eigenvector.\end{proof}
\section{The stabilizer group of a vector}
   In this section, we consider the reverse problem for the existence of eigenvectors of a Riordan matrix.
Given a nonzero vector ${\bf{h}}=(h_0,h_1,\ldots)^T$, we are interested to the Riordan matrices $A=(g,F) \in {\cal R}(\mathbb F)$ with the vector ${\bf{h}}$ as an eigenvector of $A$. Since $A{\bf{h}} = \lambda{\bf{h}}$, $\lambda\ne0$ if and only if $\big({\textstyle\frac{1}{\lambda}}\,A\big){\bf{h}} = {\bf{h}}$,
this problem is equivalent to  finding  Riordan matrices $(g,F)$
that {\it stabilize} the vector $\bf{h}$, {\it i.e.} $(g,F)\bf{h}=\bf{h}$
or  $(g,F)h(x)=h(x)$ where $h(x)$ is the generating function for the
vector ${\bf{h}}$. It is known, from the Riordan group $(g,F)$ acting on the set ${\cal F}_0$, that
for any nonzero vector $\bf{h}$, the set of all such
Riordan matrices $(g,F) \in {\cal R}(\mathbb F)$ forms a subgroup of ${\cal R}(\mathbb F)$ called the {\it
stabilizer subgroup} \cite{Tian} of ${\bf{h}}$:
$$
\stab({\bf{h}}):= \left\{ ( g,F)\in{\cal R}(\mathbb F) \mid
(g,F)\bf{h}=\bf{h}\right\}.
$$
For instance, if ${\bf{h}}=(h_0, 0, 0,  \ldots)^T$ with $h_0 \ne0$, then $\stab({\bf{h}}) = \{\,\big(1, F\big)  \mid F \in {\cal F}_1\}$, which is known as the {\it associated subgroup} \cite{Shap} of ${\cal R}(\mathbb F)$.
\vskip.3pc
By the fundamental property in (\ref{e:ftrm}), the equation
$(g,F)\bf{h}=\bf{h}$ is equivalent to $h\big(F\big) ={h\over g}$, which will be called the {\it stabilizer
equation}. Thus for any $F\in{\cal F}_1$, the Riordan matrix $\big(h/{h(F)},F\big)$ is an element of $\stab({\bf{h}})$.
One may ask that given a vector ${\bf{h}}$, is there a $F\in{\cal F}_1$ such that $(g,F)\in\stab({\bf{h}})$ for any $g\in{\cal F}_0$?
To answer this question, we will investigate the elements $(g, F) \in \stab({\bf{h}})$ from the perspective of $g\in{\cal F}_0$.
Specifically, we will delineate which $g\in {\cal F}_0$ can occur and give an explicit
formula for $F\in{\cal F}_1$ in terms of such $g$ and $h$. As a corollary, we obtain that for the $g\in {\cal F}_0$ there exists at most finitely many
$F$ such that $(g,F)\in \stab({\bf h})$.
\ms
We begin with a vector $\bf{h}$ with a single nonzero entry $h_k$ with generating function $h=h_kx^k$ with  $k\ge1$.
By the stabilizer equation it is easily shown that for any $g\in{\cal F}_0$, the Riordan matrices $(g, F) = \big(g, f_1xg^{-1/k}\big)$ with $f_1^k = 1$,  are the elements of $\stab({\bf{h}})$, because its $k^{\rm th}$ column generating function is $gF^k  = g(f_1xg^{-1/k})^k=x^k$.

Now we may assume that $\bf{h}$ has at least two nonzero entries. In the sequel, we assume that ${\bf{h}} = (h_0, 0, \ldots, 0, h_k, h_{k + 1}, \ldots)$ with $h_k\ne0$ for $k\ge1$, and we explore an element $(g,F)$ of stab($\bf h$). We note that the generating series of ${\bf h}$ can be written as
\begin{eqnarray}\label{h}
h = h_0 + \sum_{n\ge k} h_nx^n=h_0+h_kH(x)^k,\;\;h_k\ne0
\end{eqnarray}
where
$$
H(x)= x\left(1 + \frac{h_{k + 1}}{h_k}x +  \frac{h_{k + 2}}{h_k}x^2 + \cdots \right)^{1/k}\in {\cal F}_1.
$$
The compositional inverse of $H$ will be denoted by ${\overline H}$.

 \begin{lem}\label{lemma4.4} Let $g\in{\cal F}_0$ and $h = h_0 + \sum_{n\ge k} h_nx^n$, with $h_k\ne0$. If $c={h\over g}=\sum_{n\ge0}c_n$ then $c_1=\cdots=c_{k-1}=0$ if and only if $g_1=\cdots=g_{k-1}=0.$ Moreover, if $c_i=0$ or $g_i=0$ for $i=1,\ldots,k-1$ then $h_k=c_0g_k+c_kg_0$ where $c_0=h_0/g_0$.
   \end{lem}
  \begin{proof} Consider $h=cg= h_0 + \sum_{n\ge k} h_nx^n$. Then the proof is straightforward from the convolution rule:
$$
h_n=[x^n]h=[x^n]cg=\sum_{i=0}^nc_ig_{n-i}.
$$
\end{proof}

The stabilizer equation leads to the following stabilizer theorem.

\begin{thm} \label{StabThm}{\rm (The stabilizer theorem)} Let $A=(g,F)\in{\cal R}(\mathbb F)$. Assume that ${\bf h}$ is a given vector with generating
series $h=h_0+h_kH(x)^k$ with $h_k\ne0$ in (\ref{h}), and let $D=\left(h_k^{-1}({h\over g}-h_0)\right)^{1/k}$.
\begin{itemize}
\item[{\rm (a)}] If $h_0 \ne 0$, then $A\in \stab({\bf{h}})$ if and
only if $g=1+\sum_{n\ge k}g_nx^n$ and $F={\overline H}(D)$ where $f_1=[x]F$ is a root of $f_1^k =(h_k - h_0g_k)/h_k$.
\item[{\rm (b)}] If $h_0 = 0$, then $A\in \stab({\bf{h}})$ if and
only if $g=g_0+\sum_{n\ge k}g_nx^n$ and $F={\overline H}(D)$ with $h_0 = 0$ where $f_1=[x]F$ is a root of $f_1^k =1/g_0$.
\end{itemize}
\end{thm}
\begin{proof} (a) Let $h_0 \ne 0$. By the stabilizer equation, $A\in \stab({\bf{h}})$ if and
only if
\begin{eqnarray}\label{h/g}
{h\over g}=h(F)=h_0+h_kH(F)^k.
\end{eqnarray}
Now let $c=h/g=\sum_{n\ge0}c_nx^n$. Since
\begin{eqnarray}\label{c}
c=h(F)=h_0+\left(h_kf_1^k \right)x^k+{\rm(higher\; power)},
\end{eqnarray}
it implies that $c_0=h_0$, $c_1=\cdots=c_{k-1}=0$ and $c_k=h_kf_1^k \ne 0$.
Since $c_0=h_0/g_0=h_0$ we have $g_0=1$. Thus by Lemma \ref{lemma4.4}
we obtain $g=1+g_kx^k+\cdots$ where $g_k=(h_k-c_k)/h_0$. Since $c_k=h_kf_1^k$ we have $f_1^k =(h_k - h_0g_k)/h_k$.
Moreover, it follows from (\ref{h/g}) that
 $$H(F)^k={1\over h_k}\left({h\over g}-h_0\right)=D^k.$$
By the extension of Niven's Theorem in Lemma \ref{roots}, $H(F)=D\in{\cal F}_1$. Therefore, $A\in \stab({\bf{h}})$ if and only if $g=1+\sum_{n\ge k}g_nx^n$ and $F={\overline H}(D)$ where $f_1^k =(h_k - h_0g_k)/h_k$.
 \ms
(b) If $h_0 = 0$ then the result follows from the similar argument used in (a).
\end{proof}

\begin{ex}{\rm Let ${\bf{h}}=(1,1,\ldots)^T$. By (\ref{h}) we have $k=1$, $h={1\over
1-x}$, $H = \frac{x}{1-x}$, and $\overline{H} = \frac{x}{1+x}$. Since $D={h\over g}-1={1\over (1-x)g}-1$ for $g\in{\cal F}_0$ with $g_1\ne1$, it follows from (a) of Theorem \ref{StabThm} that $F=\overline{H}(D)=1-g+xg$ where $f_1=1-g_1$. Thus we obtain
 $${\rm stab}({\bf h}) =\left\{(g,F) \in {\cal R}(\F)| F=1-g+xg,\;g_1\ne1\right\},$$
which is known as the {\it stochastic subgroup} \cite{JL-Nkwanta2013}.}
\end{ex}

The following corollary is an immediate consequence of Theorem \ref{StabThm}. It asserts that for the $g\in {\cal F}_0$ there exists at most finitely many
$F$ such that $(g,F)\in \stab({\bf h})$.
\begin{cor} Let $h=h_0+\sum_{n\ge k}h_nz^n$ with $h_k\ne0$ for $k\ge1$. Given $g\in{\cal F}_0$ that satisfies (a) or (b) of Theorem \ref{StabThm},
let
$$
S_g=\left\{(g,F) \in stab({\bf{h}})\mid F\in{\cal F}_1\right\}
$$
be the subset of the stabilizer subgroup of ${\bf{h}}$. Then $|S_g|\le k$. Moreover, we have the following:
\begin{itemize}
\item[{\rm(a)}] If $h_0 \ne 0$ then $S_g$ is in one-one correspondence
with the set of $k^{\rm th}$ roots of ${h_k - h_0g_k\over h_k}$;
\item[{\rm(b)}] If $h_0=0$ then $S_g$ is in one-one correspondence
with the set of $k^{\rm th}$ roots of ${1\over g_0}$;
\item[{\rm(c)}] If $k = 1$, or if $\F = \r$ and $k$ is odd, then there exists a unique $F\in{\cal F}_1$ such that $(g,F) \in
stab({\bf{h}})$.
\end{itemize}
\end{cor}
\ms
\noindent{\bf Acknowledgement} The authors would like to thank the reviewer for  valuable comments that led to a big improvement of this paper.

\end{document}